\documentclass[11pt]{amsart}
\hoffset= - 0.75 in

\usepackage{amsfonts, amsmath, amssymb, amsthm}
\usepackage{latexsym}
\usepackage{graphicx}
\newtheorem{thm}{Theorem}[section]

\newtheorem{defn}[thm]{Definition}
\newtheorem{ex}[thm]{Example}

\newtheorem{lem}[thm]{Lemma}
\newtheorem{cor}[thm]{Corollary}
\newtheorem{prop}[thm]{Proposition}

\newcommand{\RR}{\mathbb{R}}
\newcommand{\QQ}{\mathbb{Q}}
\newcommand{\CC}{\mathbb{C}}
\newcommand{\FF}{\mathbb{F}}
\newcommand{\NN}{\mathbb{N}}
\newcommand{\M}{\mathcal{M}}
\newcommand{\C}{\mathcal{C}}
\newcommand{\TP}{\mathbb{TP}}

\newcommand{\Sym}{\mathcal{S}}
\newcommand{\bin}[2]{{#1\choose #2}}
\newcommand{\conv}{\operatorname{conv}}

\begin{document}

\title{On the rank of a tropical matrix}
\thanks{This work was conducted 
during the Discrete and Computational 
Geometry semester at M.S.R.I., Berkeley.
Mike Develin held the AIM Postdoctoral Fellowship 2003-2008 and
Bernd Sturmfels held the MSRI-Hewlett Packard
Professorship 2003/2004. 
Francisco Santos was partially supported by 
the Spanish Ministerio de Ciencia y Tecnolog\'{\i}a (grant BFM2001-1153).
Bernd Sturmfels was partially supported by the National Science Foundation
(DMS-0200729).}

\author[Develin, Santos, and Sturmfels]
{Mike Develin, Francisco Santos, and Bernd Sturmfels}
\address{Mike Develin, American Institute of Mathematics, 360 Portage
Ave., Palo Alto, CA 94306, USA}
\email{develin@post.harvard.edu}
\address{Francisco Santos, Depto. de Matem\'aticas, Estad\'{\i}stica y Computaci\'on, Universidad de Cantabria, E-39005 Santander, SPAIN}
\email{fsantos@unican.es}
\address{Bernd Sturmfels, Department of Mathematics, University of California, Berkeley, CA 94720, USA}
\email{bernd@math.berkeley.edu}
\date{\today}

\begin{abstract}
This is a foundational paper in tropical linear algebra, which is
linear algebra over the min-plus semiring. We introduce and compare three
natural definitions of the rank of a matrix, called the
Barvinok rank, the Kapranov rank and the tropical rank. We demonstrate how these notions arise 
naturally in polyhedral and algebraic geometry, and we show that they differ in general. 
Realizability of matroids plays a crucial role here.
Connections to optimization are also discussed.
\end{abstract}

\maketitle

\section{Introduction}
The rank of a matrix $M$ is one of the most important notions in
linear algebra.  This number can be defined in many different ways. In
particular, the following three definitions are equivalent:

{\em
\begin{itemize}

\item The {\em rank} of $M$ is the smallest integer $r$ for which $M$ can be written
as the sum of $r$ rank one matrices. A matrix has rank 1 if it is the product of a column vector and a 
row vector.
\item The {\em rank} of  $M$ is the smallest dimension of any linear space containing the 
columns of $M$.
\item The {\em rank} of $M$ is the largest integer $r$ such that $M$ has a
non-singular $r\times r$ minor.
\end{itemize}
}

Our objective is to examine these familiar definitions over an
algebraic structure which has no additive inverses.  We work over the
\emph{tropical semiring} $(\RR, \oplus, \odot)$ whose arithmetic
operations are
\[
 a \, \oplus \, b \,\,\, := \,\, {\rm min}(a,b) \qquad \hbox{and}
\qquad a \, \odot \, b \,\,\, := \,\, a + b.
\]
The set $\RR^d$ of real $d$-vectors and the set $\RR^{d \times n}$ of
real $d \times n$-matrices are semimodules
over the semiring $(\RR, \oplus, \odot)$. The operations of matrix addition and
matrix multiplication are well defined.
All our definitions of rank make sense over the
tropical semiring $(\RR, \oplus, \odot)$:

\begin{defn} \label{defn:brank}
The {\em Barvinok rank} of a matrix
$M\in \RR^{d\times n}$ is the smallest integer $r$ for which $M$ can be written
as the tropical sum of $r$ matrices, each 
of which is the tropical product of a $d \times 1$-matrix and a $1 \times n$-matrix. 
\end{defn}

\begin{defn}
\label{defn:krank}
 The {\em Kapranov rank} of a matrix
$M\in \RR^{d\times n}$ is the 
smallest dimension of any tropical linear space 
(to be defined in Definition \ref{troplinspace})
containing the columns of $M$.
\end{defn}
\begin{defn} \label{defn:trank1}
A square matrix $M = (m_{ij})\in \RR^{r\times r}$
is \emph{tropically singular} if the minimum in
\[
{\rm det}(M) \quad := \quad
\bigoplus_{ \sigma \in \Sym_r}  m_{1 \sigma_1} \odot
m_{2 \sigma_2} \odot \cdots \odot m_{r \sigma_r}
\quad = \quad
{\rm min} \bigl\{  \,m_{1 \sigma_1} +
m_{2 \sigma_2} + \cdots + m_{r \sigma_r} \, : \,
\sigma \in \Sym_r \,\bigr\} 
\]
is attained at least twice. Here $\Sym_r$ denotes the
symmetric group on $\{1,2,\ldots,r\}$.
The {\em tropical rank} of a matrix
$M\in \RR^{d\times n}$ is the largest integer $r$ such that $M$ has a
non-singular $r\times r$ minor.
\end{defn}

These three definitions are easily seen to agree for $r=1$, but in general they
are not equivalent:

\begin{thm} \label{thm:main}
For every matrix $M$ with entries in the tropical 
semiring $(\RR, \oplus, \odot) $, we have
\begin{equation}
\label{eqn:rankinequalities}
\hbox{tropical rank}\,(M) \quad \leq \quad
\hbox{Kapranov rank}\,(M) \quad \leq \quad
\hbox{Barvinok rank}\,(M) .
\end{equation}
Both of these inequalities can be strict.
\end{thm}

The proof of  Theorem~\ref{thm:main}  consists of
Propositions \ref{prop:KBcompare}, \ref{prop:TKcompare},
\ref{prop:matroidtrank} and Theorem \ref{prop:matroidkrank}
in this paper. As we go along,  several alternative characterizations of
the Barvinok, Kapranov and tropical ranks
will be offered. One of them arises from the fact
that every $d \times n$-matrix $M$ defines a tropically linear map
$\RR^n \rightarrow \RR^d$. 
The image of $M$ is a polyhedral complex in $\RR^d$.
Following \cite{DS}, we identify this polyhedral complex with its
image in the tropical projective space
$\TP^{d-1} = \RR^d/\RR(1,1,\ldots,1)$. This image is 
the \emph{tropical convex hull} of (the columns of) $M$ as in \cite{DS}.
Equivalently, this \emph{tropical polytope} is the set  of all tropical linear
combinations of the columns of $M$. We show in Section 4 that the tropical
rank of $M$ equals the dimension of this tropical polytope plus one, thus justifying the definition of
the vanishing of the determinant given in Definition~\ref{defn:trank1}. 

The discrepancy between Definition \ref{defn:trank1} and Definition \ref{defn:krank} comes from 
the crucial
distinction between tropical polytopes and tropical linear spaces, as explained in \cite[\S 1]{RGST}. The latter are
described in~\cite{SS} where it is shown that they are parametrized by the tropical Grassmannian.
That the two inequalities in Theorem~\ref{thm:main} can 
be strict corresponds to two  facts about tropical geometry 
which are unfamiliar from classical geometry. Strictness of the first
inequality corresponds to the fact that a 
point configuration in tropical space can have a
$d$-dimensional convex hull but not lie in any $d$-dimensional affine 
subspace. Strictness of the
second inequality corresponds to the fact that a 
point configuration in a $d$-dimensional subspace need
not lie in the convex hull of $d+1$ points.

We start out in Section 2 by studying
the Barvinok rank (Definition \ref{defn:brank}).
This notion of rank arises in the context of
combinatorial optimization~\cite{BJW,But,CRW}.
In Section 3 we study the Kapranov rank (Definition \ref{defn:krank}). This notion is the most
natural one from the point of view of algebraic geometry, where tropical arithmetic arises as the ``tropicalization'' of
arithmetic in a power series ring. It has good algebraic and geometric properties but is difficult to characterize
combinatorially; for instance, it depends on the base field of the power series ring, which here we take to be 
the complex numbers $\CC$, unless otherwise stated.

In Section 4 we study the tropical rank (Definition \ref{defn:trank1}).  This is the best notion of
rank from a geometric and combinatorial perspective. For instance, it can be expressed in terms of
regular subdivisions of products of simplices~\cite{DS}. In
Section 5, we use this characterization to show that 
the tropical
and Kapranov ranks agree when either of them is equal to ${\rm min}(d,n)$. 

Section 6 is devoted to another case where the Kapranov and tropical ranks agree, namely when either of
them equals two. The set of $d \times n$-matrices enjoying this property is the space of trees with $d$
leaves and $n$ marked points. This space is studied in the companion paper \cite{Dev}.

The second inequality of Theorem~\ref{thm:main}
is strict for many matrices (see Proposition \ref{prop:barvinokCn} for examples),
but it requires more effort to find matrices for which the first inequality is strict.  Such matrices are
constructed in Section 7 by relating Kapranov rank to realizability of matroids. 

Our definition of ``tropically non-singular'' is equivalent to what is called
``strongly regular'' in the literature on the min-plus algebra \cite{BH,C}.
The resulting notion of tropical rank, as well as the notion
of Barvinok rank, have previously appeared in this literature.
In fact, linear algebra in the tropical semiring has been 
called  ``the linear algebra of combinatorics''
by Butkovic~\cite{But}.
In the final
section of the paper we revisit some of that literature, which is concerned mainly with algorithmic
issues, and relate it to our results.  We also point out several (mostly algorithmic)
open questions.

Summing up, the three definitions of rank studied in this paper generally disagree, and they have
different flavors (combinatorial, algebraic, geometric). But they all share some of the familiar
properties of matrix rank over a field. The following properties are easily checked for each of the 
three definitions of rank: the rank of a matrix and its transpose are the
same; the rank of a minor cannot exceed that of the whole matrix; the rank is invariant under
(tropical) multiplication of rows or columns by constants, and under insertion of a row or column
obtained as a combination of others; the rank of $M \oplus N$ is at most the
sum of the ranks of $M$ and $N$; the rank of $(M\,|\,N)$ is at least the ranks of $M$ and of $N$ and at 
most the sum of their ranks; and the rank of $M \odot N$ is at most the minimum of the ranks
of $M$ and $N$. 

\section{The Barvinok rank}

The Traveling Salesman Problem can be solved in polynomial time
if the distance matrix is the tropical sum of $r$ matrices of
tropical rank one (with $\oplus$ as ``max'' instead of ``min'').
This result was proved by Barvinok, Johnson and Woeginger \cite{BJW}, 
building on earlier work of Barvinok.
This motivates our definition of {\em Barvinok rank} as the smallest $r$ for which $M\in 
\RR^{d\times n}$
is expressible in this fashion. Since matrices of tropical rank one are of the form $X\odot Y^T$,
for two column vectors $X \in \RR^d$ and $Y \in \RR^n$, this is equivalent to saying that $M$ has a representation
\begin{equation}
\label{eqn:brank}
 M \quad = \quad
X_1 \odot Y_1^T \,\,\oplus \,\,
X_2 \odot Y_2^T \,\,\oplus \,\, \cdots \,\,\oplus \,\,
X_r \odot Y_r^T.
\end{equation}

\begin{figure}
\includegraphics{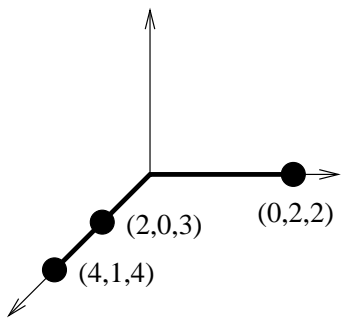}
\caption{\label{fig:rk2}
A tropical line in $\TP^2$, and a one-dimensional tropical polytope.}
\end{figure}

For example, the following equation shows a
$3 \times 3$-matrix which has Barvinok rank two:
$$
\begin{pmatrix}
0 & 4 & 2 \\
2 & 1 & 0 \\
2 & 4 & 3
\end{pmatrix} \quad = \quad
\begin{pmatrix} 
0 \\
2 \\
2
\end{pmatrix}
\odot (0,4,2) \,\,\,\, \oplus \,\,\,\,
\begin{pmatrix}
3 \\
0 \\
3
\end{pmatrix} 
\odot (2,1,0).
$$
This matrix also has tropical rank $2$ and
Kapranov rank $2$ because the matrix
is tropically singular. The column vectors lie on the tropical line
in $\TP^2 = \RR^3/\RR(1,1,1)$ defined by
$\,2\odot x_1 \oplus 3\odot x_2 \oplus 0\odot x_3 $, depicted in Figure~\ref{fig:rk2}. Their convex hull, 
darkened, is a subset of the line and thus one-dimensional.

We next present two reformulations of the definition 
of Barvinok rank:  in terms of tropical convex hulls
as introduced in  \cite{DS}, and via a ``tropical morphism'' 
between matrix spaces.

\begin{prop}
\label{prop:barvinok}
Let $M$ be a real $d \times n$-matrix. The following
properties are equivalent:
\begin{itemize}
\item[(a)] $M$ has Barvinok rank at most $r$.
\item[(b)] The columns of $M$
lie in the tropical convex hull of $r$ points in $\TP^{d-1}$.
\item[(c)] There are matrices $X\in \RR^{d\times r}$ and $Y\in
\RR^{r\times n}$ such that $M=X\odot Y$.
 Equivalently, $M$ lies in the image of the following
tropical morphism, which is defined by matrix multiplication:
\begin{equation}
\label{eqn:morphism}
\phi_r\,\,:\,\, \RR^{d \times r}  \times  \RR^{r \times n}  \, \rightarrow \,
\RR^{d \times n} \,\, , \,\,\,
(X,Y) \,\mapsto \,X \odot Y .
\end{equation}
\end{itemize}
\end{prop}

\begin{proof}
Let $M_1,\dots, M_n \in \RR^{d}$ be the column vectors of $M$. Let
$X_1,\dots,X_r\in \RR^d$ and $Y_1,\dots,Y_r\in\RR^n$ be the columns of
two unspecified matrices $X\in \RR^{d\times r}$ and $Y\in \RR^{n\times
r}$. Let $Y_{ij}$ denote the $j$th coordinate of $Y_i$.
The following three algebraic identities are easily seen to be equivalent:
\begin{itemize}
\item[(a)] 
$M \, = \,
X_1 \odot Y_1^T \,\oplus \,
X_2 \odot Y_2^T \,\oplus \, \cdots \,\oplus \,
X_r \odot Y_r^T,$
\item[(b)] $M_j \, = \,
Y_{1j}\odot X_1 \,\oplus \,
Y_{2j} \odot X_2 \,\oplus \, \cdots \,\oplus \,
Y_{rj} \odot X_r $ for all $j = 1,\ldots,n$, and
\item[(c)] $M=X\odot Y^T$.
\end{itemize}
Statement (b) says that each column vector of $M$ lies in the tropical convex
hull of $X_1,\ldots,X_r$. The entries of the matrix $Y$
are the multipliers in that tropical linear combination.
This shows that the three conditions (a), (b) and (c)
in the statement of the proposition are equivalent.
\end{proof}

Part (b) of Proposition \ref{prop:barvinok} suggests that the Barvinok
rank of a tropical matrix is more an analogue of the non-negative rank
of a matrix than of the usual rank.  Recall (e.g.~from \cite{CR}) that
the \emph{non-negative rank} of a real non-negative matrix $M\in
\RR^{d\times n}$ is the smallest $r$ for which $M$ can be written as a
product of non-negative matrices of format $d\times r$ and $r\times
n$.  Equivalently, it is the smallest $r$ for which the columns (or
rows) of $M$ lie in the positive hull of $r$ non-negative
vectors. Compare this with the formulation of Barvinok rank given in
Proposition~\ref{prop:barvinok} (b); this closer connection comes from
the fact that tropical
linear combinations yield an object
more analogous to a ``positive span'' or ``convex hull''~\cite{DS,
RGST} than a linear span.  For more information on non-negative rank
see \cite{CR}, and for the connection to rank over other semigroup
rings see \cite{GP}.

By Proposition \ref{prop:barvinok}, 
the set of all Barvinok matrices  of rank $\leq r$
is the image of the tropical morphism $\phi_r$. In particular,
this set is a polyhedral fan in $\RR^{d \times n}$. This fan has 
interesting combinatorial structure, 
even for $r=2$. These fans are
discussed in more detail in~\cite{Ard} and~\cite{Dev}.


We next present an example of
a matrix which shows that the
Barvinok rank can be much larger than the
other two notions of rank. The matrix
to be considered is the \emph{classical identity matrix} 
\begin{equation}
\label{eqn:LooksLikeUnitMatrix}
C_n \quad = \quad
\begin{pmatrix}
1 & 0 & 0 & \cdots & 0 \\
0 & 1 & 0 & \cdots & 0 \\
0 & 0 & 1 & \cdots & 0 \\
\vdots &\vdots &\vdots & \ddots & \vdots \\
0 & 0 & 0 & \cdots & 1
\end{pmatrix}.
\end{equation}
This looks like the unit matrix
(in classical arithmetic) but it is
far from being a unit matrix in
tropical arithmetic, where
$0$ is the neutral element for $\odot$ and $\infty$ is
the neutral element for $\oplus$. After obtaining the
following result, we learned that the  same
calculation had already been done in~\cite{CRW}.

\begin{prop}
\label{prop:barvinokCn}
The Barvinok rank of the matrix $C_n$ is
 the smallest integer $r$ such that
\[
 n \, \le \, \bin{r}{\left\lfloor\frac{r}{2}\right\rfloor}.
\]
\end{prop}

\begin{proof}
Let $r$ be an integer
and assume that $n\le {r \choose \lfloor r/2\rfloor}$.
We first show that $\,\hbox{Barvinok~rank$\,(C_n)$} \leq r$.
Let $S_1,\dots,S_n$ be distinct subsets of $\{1,\dots,r\}$
each having cardinality $\lfloor r/2\rfloor$.
For each $k\in 1,\dots, r$, we define an
$n \times n$-matrix $\,X_k = (x_{ij}^k)\,$
with entries in  $\{0,1,2\}$ as follows:
$$\hbox{ $x_{ij}^k=0\,\,$ if $\,\, k\in S_i \backslash
S_j \,$, $\quad
x_{ij}^k=2\,\,$ if $\,\,k\in S_j\backslash S_i$, \ and 
$\quad x_{ij}^k=1\,$ otherwise}.
$$
The matrix $X_k$ has tropical rank one. To see this, let
 $\,V_k\in \{0,1\}^n\, $ denote the vector with $i$th
coordinate equal to one or zero depending on whether $k$
is an element of $ S_i$ or not.
Then we have
\[
X_k \quad = \quad V_k^T \odot (\,{1}\odot(-V_k)\,).
\]
To prove $\,\hbox{Barvinok~rank$\,(C_n)$} \leq r$, it now suffices to 
establish the identity
\[
C_n \quad = \quad X_1 \,\oplus \,X_2 \,\oplus\, \cdots \,\oplus \, X_{r}.
\]
Indeed, all diagonal entries of the matrices
on the right hand side are $1$, and the off-diagonal
entries  (for $i \not= j$) of the right hand side are
$\,{\rm min}(x^1_{ij}, x^2_{ij}, \ldots, x^r_{ij}) \, = \, 0$,
because $\, S_i \backslash S_j\,$ is non-empty.

To prove the converse direction, we consider an arbitrary representation
\[
C_n \quad = \quad Y_1\,\oplus\, Y_2 \,\oplus \,\cdots\,\oplus \,Y_{r}
\]
where the matrices $Y_k = (y^k_{ij})$ have tropical rank one.
For each $k$ we set $\,T_k \,:= \,\{(i,j) : y^k_{ij}=0\}$.
Since the matrices $Y_k$ are non-negative and have
tropical rank one, it follows that each
 $T_k$ is a product $I_k\times J_k$,
where $I_k$ and $J_k$ are subsets of $\{1,\ldots,n\}$.
Moreover, we have
$\,I_k\cap J_k=\emptyset\, $ because 
 the diagonal
entries of $Y_k$ are not zero.
For each $i=1,\dots, n$ we set
\[
S_i \,:= \,\{ k : i\in I_k\} \quad \subseteq \,\,\,\{1,\dots,r\}.
\]
We claim that no two of the sets $S_1,\dots,S_n$ are contained in one
another. Sperner's Theorem \cite{AZ} then proves that
$n\le {r\choose \lfloor r/2\rfloor}$. To prove the claim, observe that
if $S_i\subset S_j$ then the entry $y^k_{i,j}$ cannot be zero for any
$k$. Indeed, if $k\in S_i\subseteq S_j$ then $j\in I_k$
implies $j\not\in J_k$. And if $k\not\in S_i$ then $i\not\in I_k$.
\end{proof}

For example, $C_6$ has Barvinok rank $4$, as the following
decomposition shows:
\[
C_6=
\begin{pmatrix}
1 & 1 & 1 & 2 & 2 & 2 \\
1 & 1 & 1 & 2 & 2 & 2 \\
1 & 1 & 1 & 2 & 2 & 2 \\
0 & 0 & 0 & 1 & 1 & 1 \\
0 & 0 & 0 & 1 & 1 & 1 \\
0 & 0 & 0 & 1 & 1 & 1 \\
\end{pmatrix}
\oplus
\begin{pmatrix}
1 & 1 & 0 & 0 & 0 & 1 \\
1 & 1 & 0 & 0 & 0 & 1 \\
2 & 2 & 1 & 1 & 1 & 2 \\
2 & 2 & 1 & 1 & 1 & 2 \\
2 & 2 & 1 & 1 & 1 & 2 \\
1 & 1 & 0 & 0 & 0 & 1 \\
\end{pmatrix}
\oplus
\begin{pmatrix}
1 & 0 & 1 & 0 & 1 & 0 \\
2 & 1 & 2 & 1 & 2 & 1 \\
1 & 0 & 1 & 0 & 1 & 0 \\
2 & 1 & 2 & 1 & 2 & 1 \\
1 & 0 & 1 & 0 & 1 & 0 \\
2 & 1 & 2 & 1 & 2 & 1 \\
\end{pmatrix}
\oplus
\begin{pmatrix}
1 & 2 & 2 & 2 & 1 & 1 \\
0 & 1 & 1 & 1 & 0 & 0 \\
0 & 1 & 1 & 1 & 0 & 0 \\
0 & 1 & 1 & 1 & 0 & 0 \\
1 & 2 & 2 & 2 & 1 & 1 \\
1 & 2 & 2 & 2 & 1 & 1 \\
\end{pmatrix}.
\]
Similarly, $C_{36}$ has Barvinok rank $8$, even though all its 
$35\times 35$ minors have Barvinok rank $7$ (and its $8\times 8$
minors have Barvinok rank at most $5$). Asymptotically, 
$$ \hbox{Barvinok~rank$\,(C_n)\quad \sim \quad \log_2 n$}. $$
We will see in Examples \ref{ex:Krank2} and \ref{ex:trank} that the
Kapranov rank and tropical rank of
$C_n$ are both two.

\section{The Kapranov rank}

The tropical semiring has a strong connection to power series rings
and their algebraic geometry. We review the basic setup
from \cite{SS, Stu}. Let $K= \CC\{\! \{t\} \! \}$ be the field of Puiseux 
series with complex coefficients. The elements in $K$ are
formal power series
$\, f \, = \, c_1 t^{a_1} + c_2 t^{a_2} + \cdots$, where
$\,a_1 < a_2 < \cdots \,$ are rational numbers that have
 a common denominator. 
Let ${\rm deg}: K^* \rightarrow \QQ$ be the natural valuation
sending a non-zero Puiseux series $f$ to its degree $a_1$. 
For any two elements $f,g \in K$, we have
$\,{\rm deg}(fg) \,=\,
{\rm deg}(f) + {\rm deg}(g) \,=\, {\rm deg}(f)\odot {\rm deg}(g)$.
In general we also have
${\rm deg}(f+g) = {\rm min}\,({\rm deg}(f),{\rm deg}(g)) = {\rm
deg}(f)\oplus {\rm deg}(g)$, unless there is a cancellation of leading terms.
 Thus the tropical arithmetic is naturally induced from
ordinary arithmetic in  power series fields.

The field  $K= \CC\{\! \{t\} \! \}$ is algebraically closed
of characteristic zero. 
If $I$ is any ideal in $K[x_1,\ldots,x_d]$
then we write $V(I)$ for its variety in the
$d$-dimensional algebraic torus $(K^*)^d$.
Thus the elements of $V(I)$ are vectors
$x(t) = (x_1(t),\ldots,x_d(t))$ where each
$x_i(t)$ is a Puiseux series and 
$\, f(x(t)) = 0 \,$ for each polynomial $f \in I$.
Let us now enlarge the field $K$ and allow all formal power series
$\, f \, = \, c_1 t^{a_1} + c_2 t^{a_2} + \cdots$ where
the $a_i$ can be real numbers, not just rationals.
We denote this larger field by $\tilde K$ and we write
$\tilde V(I)$ for the variety in $(\tilde K^*)^d$ defined by $I$.
The degree map can be applied coordinatewise,
giving rise to a map which takes vectors
of non-zero power series into $\RR^d$:
$$  {\rm deg} \, : \, (\tilde K^*)^d \,\,\, \rightarrow \,\,\, \RR^d \,,\,\,
(f_1(t), \ldots,f_d(t)) \,\mapsto \, \bigl({\rm deg}(f_1),\ldots,
{\rm deg}(f_d) \bigr).   $$

We define the \emph{tropical variety} of $I$,
denoted $\mathcal{T}(I) \subset \RR^d$, to be the
image of $\tilde V(I)$ under the map ${\rm deg}$.
In \cite{SS, Stu},  the following alternative description of
the tropical variety is given:

\begin{thm} \label{thm:twice}
The tropical variety $\mathcal{T}(I)$ is the set of
vectors $w \in \RR^n$ such that the initial
ideal $ {\rm in}_w(I)  = \langle  {\rm in}_w(f) : f \in I \rangle$
 contains no monomial. The dimension of $\mathcal{T}(I)$ is
the (topological) dimension of $V(I)$.
\end{thm}

The first statement in Theorem \ref{thm:twice} is due to
Misha Kapranov (in the special case when $I$ is a principal ideal) 
and the third author (for arbitrary ideals $I$, in \cite{Stu}).
A complete proof can be found in \cite{SS}.
The second statement  in Theorem \ref{thm:twice} 
is due to Bieri and Groves \cite{BGr}.
An elementary proof of this result, and the fact that 
$\mathcal{T}(I)$ is a polyhedral fan, appears in \cite[\S 9]{Stu}.

We defined Kapranov rank to be the smallest dimension 
of any tropical linear space containing the
columns of $M$; now, we can make this precise by 
defining tropical linear spaces.

\begin{defn} \label{troplinspace}
A \emph{tropical linear space} in $\RR^d$ is any subset $\mathcal{T}(I)$ where
$I$ is an ideal generated by affine-linear forms
$\, a_1 x_1 + \cdots + a_d x_d + b \,$ in $\,\tilde K[x] = \tilde K[x_1,\ldots,x_d]$. Its 
\emph{dimension} is its topological dimension, which is equal to $d$ minus 
the number of minimal generators of  $I$.
\end{defn}

Note that here the scalars $a_1,\ldots,a_n, b$ are power series in $t$
with complex coefficients, the choice of the complex numbers being crucial.
If $I$ is the principal ideal  generated by one affine-linear form
$\, a_1 x_1 + \cdots + a_n x_n + b $, then $\mathcal{T}(I)$
is a \emph{tropical hypersurface}. Tropical linear spaces
were studied in \cite{SS}, where it was shown that they
are parametrized by the \emph{tropical Grassmannian}.
Every tropical linear space $L$ is a finite intersection 
of tropical hyperplanes, but not conversely, and the number
of tropical hyperplanes needed is
generally larger than the codimension of $L$.

Recall from Definition \ref{defn:krank} that the Kapranov rank of a matrix
$M\subset \RR^{d\times n}$ is the 
smallest dimension of any tropical linear space containing the columns of $M$.
It is not completely apparent in this definition that the Kapranov rank of a matrix 
and its transpose are the same,
but this follows from our next result.  Let $J_r$ denote the ideal generated by all the $(r+1) \times
(r+1)$-subdeterminants of a $d \times n$-matrix of indeterminates $(x_{ij})$. This is a prime ideal of dimension $rd
+ rn - r^2$, and the generating determinants form a Gr\"obner basis. The variety $V(J_r)$ consists of all $d \times
n$-matrices with entries in $K^*$ whose (classical) rank is at most $r$.

\begin{thm}\label{thm:kaprlift}
For a real matrix $M = (m_{ij}) \in \RR^{d \times n}$ the following statements are equivalent:
\begin{enumerate}
\item[(a)]  The Kapranov rank of $M$ is at most $r$.
\item[(b)] The matrix  $M$ lies in the tropical determinantal variety $\mathcal{T}(J_{r})$.
\item[(c)]  There exists a $d\times n$-matrix $F = \bigl(f_{ij}(t)\bigr)$ with non-zero entries
in the field $\tilde K$ such that the rank of $F$ is less than or equal to 
$r$ and
 $\,{\rm deg}(f_{ij}) = m_{ij}\,$ for all $i$ and $j$.
\end{enumerate}
\end{thm}

The  power series matrix $F$ in part (c) is called a \emph{lift} of $M$. We abbreviate
this as $\,{\rm deg}(F) = M$.

\begin{proof}
The equivalence of (b) and (c) is simply our definition of tropical variety applied to the 
ideal $J_{r}$ since, over the 
field $\tilde K$, lying in the variety of the determinantal ideal $J_{r}$ is equivalent to having rank at most 
$r$. To see that (c) implies (a), consider the linear subspace
of $\tilde K^d$ spanned by the columns of $F$. This is an 
$r$-dimensional linear space over a field, so it is defined by an ideal $I$
generated by $d-r$ linearly independent linear forms in $\tilde K[x_1,\ldots,x_d]$.
The tropical linear space $\mathcal{T}(I)$ contains all the
column vectors of $M = {\rm deg}(F)$.

Conversely, suppose that (a) holds, and let $L$ be a tropical
linear space of dimension $r$ containing the columns of $M$. 
Pick a linear ideal $I $ in  $\tilde K[x_1,\ldots,x_d]$ such 
that $L = \mathcal{T}(I)$. By applying the definition of tropical variety
to the ideal $I$, we see that each column vector of $M$
has a preimage in $\tilde V(I) \subset (\tilde K^*)^d$ under the degree map.
Let $F$ be the $d \times n$-matrix over $\tilde K$ whose columns
are these preimages. Then the column space of $F$ is contained in the
variety defined by $I$, so we have $\,{\rm rank}(F) \le r$, and
$\, {\rm deg}(F) = M$ as desired.
\end{proof}

\begin{cor}
\label{cor:kaprlift}
The Kapranov rank of a matrix $M \in \RR^{d \times n}$ is the smallest rank of any lift
of $M$.
\end{cor}

The ideal $J_1$ is generated by the $2 \times 2$-minors
$\, x_{ij} x_{kl } - x_{il} x_{kj} \,$
of the $d \times n$-matrix $(x_{ij})$.  Therefore, 
a matrix of Kapranov rank one must certainly satisfy
the linear equations $\,m_{ij} + m_{kl} = m_{il} + m_{kj}$.
This happens if and only if there exist real vectors
$\, X = (x_1,\ldots,x_d)\,$ and $\, Y = (y_1,\ldots,y_n)\,$ with
$$ m_{ij} \, = \, x_i + y_j \,\,\,\hbox{for all} \,\, i,j \quad \iff \quad
   m_{ij} \, = \, x_i \odot y_j \,\,\,\hbox{for all} \,\, i,j \quad \iff \quad
 M \, = \, X^T \odot Y . $$
Conversely, if such $X$ and $Y$ exist, we can lift $M$ to a matrix of
rank one by substituting $t^{m_{ij}}$ for $m_{ij}$. Therefore,
a matrix $M$ has Kapranov rank one if and only if it has Barvinok rank
one. In general, the Kapranov rank can be much smaller than the
Barvinok rank, as the following example shows.

\begin{ex} \label{ex:Krank2} \rm
Let $n \geq 3$ and
consider the classical identity matrix $C_n$.
It does not have Kapranov rank one, so it has 
Kapranov rank at least  two. 
Let $a_3,a_4,\dots,a_{n}$ be distinct nonzero complex numbers. Consider the matrix
\[ F_n \quad = \quad
\begin{pmatrix} 
t & 1 & t+ a_3 & t+ a_4 & \cdots & t+ a_n \\
1& t & 1+a_3t  & 1+a_4t  & \cdots & 1+a_nt  \\
t-a_3 & 1 & t & t-a_3+a_4 & \cdots & t-a_3+a_n \\ 
t-a_4  & 1 & t-a_4+a_3 &t &\cdots & t-a_4+a_n \\
\vdots & \vdots & \vdots & \vdots &\ddots & \vdots \\ 
t-a_n & 1 &  t-a_n+a_3 & t-a_n+a_4  & \cdots & t \\
\end{pmatrix}.
\]
The matrix $F_n$ has rank $2$ because the $i$-th column
(for $i \geq 3$) equals the first column plus $a_i$ times the second column.
Since $\, {\rm deg}(F_n) = C_n$, we conclude that $C_n$
has Kapranov rank two.

The two-dimensional tropical plane containing the columns of $C_n$ is the two-dimensional fan $L$ in $\RR^n$ which
consists of the $n$ cones $\{ x_i \geq x_1 = \cdots = x_{i-1} = x_{i+1} = \cdots = x_n \}$; this is the tropical variety
defined by the ideal in $K[x_1,\ldots,x_n]$ generated by $n-2$ linear forms with generic coefficients in $\CC$. Its
image in $\TP^{n-1}$ is the line all of whose tropical Pl\"ucker coordinates are zero \cite{SS}.
\end{ex}

The following proposition establishes  half of Theorem \ref{thm:main}.

\begin{prop}
\label{prop:KBcompare}
Every matrix $M \in \RR^{d \times n}$ satisfies
$\, \hbox{Kapranov rank}\,(M) \,\leq \,
\hbox{Barvinok rank}\,(M) $, and this inequality can be strict.
\end{prop}

\begin{proof}
Suppose that $M$ has Barvinok rank $r$. Write
$M = M_1 \oplus \cdots  \oplus M_r$ where each $M_i$ has
Barvinok rank one. Then $M_i$ has Kapranov rank one,
so there exists a rank one matrix $F_i$ over $\tilde K$
such that $\,{\rm deg}(F_i) = M_i$. Moreover, by multiplying the matrices $F_i$ 
by random complex numbers, we can choose $F_i$ such that
there is no cancellation of leading terms in $t$ when we 
form the matrix $\, F \, = \, F_1 +\cdots + F_r$. This means
$\, {\rm deg}(F) = M$. Clearly, the matrix $F$ has rank $\leq r$.
Theorem  \ref{thm:kaprlift} implies that $M$ has Kapranov rank $\leq r$.
Example \ref{ex:Krank2} shows that the inequality can be strict.
\end{proof}

A general algorithm for computing the Kapranov rank of a matrix $M$
involves computing a Gr\"obner basis of the determinantal ideal $J_r$. Suppose
we wish to decide whether a given real
$d \times n$-matrix $M = (m_{ij})$ has Kapranov rank $> r$. 
To decide this question, we fix any term order $\prec_M$ on
the polynomial ring $\CC[x_{ij}]$ which refines the partial ordering
on monomials given assigning weight $m_{ij}$ to the variable $x_{ij}$,
and we compute the reduced Gr\"obner basis  $\mathcal{G}$ of $J_r$ in the
term order $\prec_M$. For each polynomial $g$ in $\mathcal{G}$,
we consider its leading form ${\rm in}_M(g)$ with respect to the
partial ordering coming from $M$.  
As noted in \cite[\S 1]{OldStu}, we have
$\,{\rm in}_{\prec_M} \bigl( {\rm in}_M(g) \bigr) \, = \,
{\rm in}_{\prec_M} (g)$ for all $g \in \mathcal{G}$.

The ideal  generated by the set of leading forms 
$\, \bigl\{{\rm in}_M(g) \, : \, g \in \mathcal{G} \bigr\} \,$
is the initial ideal ${\rm in}_M(J_r)$. 
Let $x^{all}$ denote the product of all 
$dn$ unknowns $x_{ij}$. The second step in our
algorithm is to compute the saturation of the initial ideal with respect
to the coordinate hyperplanes:
\begin{equation}
\label{colonideal}
 \bigl( \, {\rm in}_M(J_r) \, \, : \,\,\langle x^{all} \rangle^\infty \,\bigr) \quad
= \quad \bigl\{ \, \, f \in \CC[x_{ij}] \, \,: \,\,
f \cdot (x^{all})^s \, \in \, J_r \,\,\, \hbox{for some} \,\, s \in \NN \bigr\}.
\end{equation}
Computing such an ideal quotient, given the generators
$\,{\rm in}_M(g) $, is a standard operation in computational
commutative algebra. It is a built-in command  in
software systems such as 
{\tt CoCoA}~\cite{cocoa}, {\tt Macaulay 2}~\cite{m2} or {\tt Singular}~\cite{sing}. The following 
is a direct consequence of Theorems \ref{thm:twice} and \ref{thm:kaprlift}.

\begin{cor}
The matrix $M$ has
Kapranov rank $> r$ if and only if (\ref{colonideal}) 
is the unit ideal $\,\langle 1 \rangle $.
\end{cor}

In view of this, the (combinatorial) 
Theorem \ref{thm:maxrank}, Theorem \ref{thm:rank2} 
and Corollary \ref{cor:nonrepresentable}
have the following commutative algebra
implications. Recall from \cite{RGST}
that a finite generating set $S$ of an ideal $I$ is a \emph{tropical basis} if,
for every weight vector $w \in \RR^n$ for which the initial ideal 
${\rm in}_w(I)$ contains a monomial, there is an $f\in S$ such that
 ${\rm in}_w(f)$ is a monomial. Every ideal $I$ in $K[x_1,\ldots,x_n]$
has a tropical basis but tropical bases are often much larger than
minimal generating sets. 

\begin{cor}
The $3 \times 3$-minors of
a matrix of indeterminates form a
tropical basis. The same holds
for the maximal minors of a matrix, but it
does not hold for the $4 \times 4$-minors
of a $7 \times 7$-matrix.
\end{cor}

We have defined Kapranov rank in terms of power series arithmetic over the 
complex field $\CC$, which is a canonical choice for
doing algebraic geometry.  However, the same definition works over any field $k$.
One can consider the Puiseux series field 
$ K=k\{ \! \{t\} \! \}$ 
with either rational or real exponents. Note that the former
is not algebraically closed if $k$ is algebraically closed
of characteristic $p$, but this need not concern us.
We denote the latter by $\tilde K$ as before.
All we need is the degree map $\, ({\tilde K}^*)^d \rightarrow \RR^d $. We make the following analogous 
definitions.

\begin{defn}
Let $ K=k\{ \! \{t\} \! \}$. A \emph{tropical linear space over $k$}
is the image under ``deg'' of any linear subspace
of the $\tilde K$-vector space $\tilde K^d$. Its \emph{dimension} is equal to the dimension of that 
linear subspace. The \emph{Kapranov rank over $k$} of a matrix $M\in \RR^{n\times d}$ is the smallest 
dimension of a tropical linear space containing the columns of $M$.
\end{defn}

Unless otherwise stated, we will concern ourselves only with Kapranov rank over the complex numbers. In
the general setting, Theorem~\ref{thm:kaprlift} is true over all fields, but
Proposition~\ref{prop:KBcompare} is true only over infinite fields because in its proof we needed to
take random coefficients. Indeed, Example~\ref{ex:finitefield} in Section~\ref{sec:rank2} shows a
matrix whose Kapranov rank over the 2-element field $\FF_2$ is greater than the Barvinok rank.  Even
over algebraically closed fields, the 
Kapranov rank of a matrix may depend on the characteristic of the field.
We will discuss this further and give examples in 
Section~\ref{sec:matroids}.

\section{The tropical rank }\label{sec:troprank}

We begin by proving the first inequality in Theorem \ref{thm:main}.
To complete the proof of  Theorem \ref{thm:main}, it remains
to be seen that the inequality can be strict.
This will be done in
in Section \ref{sec:matroids}.

\begin{prop}
\label{prop:TKcompare}
Every matrix $M \in \RR^{d \times n}$ satisfies
$\, \hbox{tropical rank}\,(M) \,\leq \,
\hbox{Kapranov rank}\,(M) $.
\end{prop}

\begin{proof}
If the matrix 
$M$ has a tropically non-singular $r\times r$ minor, then any lift of
$M$ to the power series field $\tilde K$
must have the corresponding $r\times r$-minor 
non-singular over $\tilde K$, 
since the leading exponent of its determinant occurs only
once in the sum. Consequently, no lift of $M$ to $\tilde K$ can have rank
less than $r$. By Theorem~\ref{thm:kaprlift}, this means that
the Kapranov rank of $M$ must be at least $r$.
\end{proof}

The set of all tropical linear combinations of a set of $n$ vectors in
$\RR^d$ is a polyhedral complex. It has a $1$-dimensional lineality
space, spanned by the vector $(1,\ldots,1)$, but upon quotienting out
by this $1$-dimensional space, we get a bounded subset in {\em
tropical projective space} $\TP^{d-1} = \RR^d / \RR(1,\ldots,1)$. This
set is the {\em tropical convex hull} of the $n$ given points in
$\TP^{d-1} $, and it was investigated in depth in~\cite{DS}. We review
some relevant definitions and facts.

We fix a subset $V=\{v_1,\dots,v_n\}\subseteq\RR^d$. Given a point $x\in \RR^d$, its {\em type} is the
$d$-tuple of sets $S=(S_1,\ldots,S_d)$, where each $S_j\subset \{1,\ldots,n\}$ and $i\in S_j$ if
$x_j-v_{ij}\ge x_k-v_{ik}$ for all $k\in \{1,\ldots,n\}$. Let $X_S$ be the region consisting of points
with type $S$; then according to~\cite[Theorem 15]{DS}, the tropical convex hull of $V$ equals the
union of the bounded regions $X_S$, which are precisely those regions for which each $S_j$ is nonempty. (If $x$ is a point in the 
tropical convex hull with type $S$, then expressing $x$ as a linear combination of the $v_i$'s, we have $i\in S_j$ if the 
contribution of $v_i$ is responsible for the $j$-th coordinate of $x$.) 
Indeed, (the topological closures of) these regions provide a polytopal 
decomposition of the tropical convex hull of $V$. Note that 
by definition, any type has the property that each $i\in \{1,\ldots,n\}$ is in some $S_j$.

The dimension of a particular cell $X_S$ of the tropical polytope 
can be easily computed from the combinatorics of the $d$-tuple $S$: 
let $G_S$ be the graph which has vertex set ${1,\ldots,d}$, with $i$ and
$j$ connected by an edge if $S_i\cap S_j$ is nonempty. The 
dimension of $X_S$ is one less than the
number of connected components of the graph $G_S$.

Recall from Definition \ref{defn:trank1} that the
tropical rank of a matrix is the size of the 
largest non-singular square minor, and that an 
$r\times r$ matrix $M$ is non-singular if 
 $\, \bigodot_{i=1}^r  M_{\sigma(i),i} \, = \,
\sum_{i=1}^r M_{\sigma(i),i}\,$
 achieves its minimum only once 
as $\sigma$ ranges over the symmetric group $\Sym_r$. 
Here is another characterization.

\begin{thm}
\label{thm:detrank}
Let $M\subset \RR^{d\times n}$ be a matrix. Then the tropical rank of 
$M$ is equal to one plus the 
dimension of the tropical  convex hull of the columns of $M$,
viewed as a  tropical polytope in $\TP^{d-1}$.
\end{thm}

\begin{proof}
Let $V=\{v_1,\ldots,v_n\}$ be the set of columns of $M$, and let $P
= {\rm tconv}(V)$ be their  tropical convex hull in $\TP^{d-1}$.
Suppose that $r$ is the tropical rank of $M$, that is,
there exists a tropically non-singular $r \times r$-submatrix
of $M$, but all larger square submatrices are tropically singular.

We first show that $\, {\rm dim}(P) \geq r-1$. We fix a
non-singular $r \times r$-submatrix $M'$ of $M$. Deleting the
rows outside $M'$ means projecting $P$ into $\TP^{r-1}$,
and deleting the columns outside $M'$ means passing to a 
tropical subpolytope $P'$ of the image. Both operations
can only decrease the dimension, so it suffices to show
${\rm dim}(P') \geq r-1$. Hence, we can assume that $M$ is itself a
tropically non-singular $r \times r$-matrix. Also, without loss
of generality, we can assume that the minimum over $\sigma\in \Sym_r$ of
\begin{equation}\label{eqn:detsingular}
f(\sigma) \quad = \quad \sum_{i=1}^r v_{i,\sigma(i)}
\end{equation}
is uniquely achieved when $\sigma$ is the identity element $e\in \Sym_r$.
We now claim that the cell $X_{(\{1\},\ldots,\{r\})}$ exists; to do this, we need to demonstrate that 
there exists a point with type $(\{1\}, \ldots, \{r\})$. 

The inequalities which must be valid on this cell are  $\,x_k-x_j
\le v_{jk}-v_{jj}\,$ for $j\neq k$. We claim that these inequalities define a full-dimensional region.
Suppose not; then, by Farkas' Lemma, there exists
a non-negative linear combination of the inequalities
 $\,x_k-x_j \le v_{jk}-v_{jj}\,$ which equals
$0\le c$ for some non-positive real $c$. This
linear combination would imply that some other $\sigma\in
\Sym_r$ has $f(\sigma)\le f(e)$, a contradiction. So this cell is full-dimensional; it follows 
immediately that picking a point in its interior yields a point with type $(\{1\}, \ldots, \{r\})$, 
since because these inequalities are all strict, no other type-inducing inequalities can hold. 

For the converse, suppose that $\, {\rm dim}(P) \geq r$.
Pick a region $X_S$ of dimension $r$, and assume by translating the points (which adds a constant to 
each row of $X_S$, not changing the rank of the matrix) that $(0,\ldots,0)$ is in 
$X_S$, so that the only inequalities valid on 0 are those given by $S$.
The graph $G_S$ has $r+1$ 
connected
components, so we can pick $r+1$ elements of
$\{1,\ldots,n\}$ of which no two appear in a common $S_j$. 
Assume without loss of generality that this set 
is $\{1,\ldots,r+1\}$, and again without loss of generality rearrange the labeling of the coordinates 
so that $i \in S_j$ if and only if $i=j$,
for $1 \leq i,j \leq r+1$.

We now claim that the square
submatrix consisting of the first $r+1$ rows and columns of $M$ is
tropically non-singular. Indeed, we have (using the definition of $f(\sigma) $
given in~$(\ref{eqn:detsingular})$):
$$
f(\sigma)-f(e) \quad = \quad
\sum_{i=1}^{r+1} v_{i,\sigma(i)}-\sum_{i=1}^{r+1} v_{ii}
\quad = \quad \sum_{i=1}^{r+1} \bigl( v_{i,\sigma(i)}-v_{ii} \bigr), $$
but whenever $\sigma(i)\neq i$, $v_{i,\sigma(i)}-v_{ii}>0$ since $i\in
S_i$ and $i\notin S_{\sigma(i)}$ for the point 0. Therefore, if $\sigma$ is not the
identity, we have $f(\sigma)-f(e)>0$, and $e$ is the unique
permutation in $\Sym_{r+1}$ minimizing the expression (\ref{eqn:detsingular}). So
$M$ has tropical rank at least $r+1$.
This is a contradiction, and we conclude
 $\, {\rm dim}(P) = r-1$.
\end{proof}

We next present a combinatorial formula
for the tropical rank of a  zero-one matrix, or 
any matrix which has only two distinct entries.
We define
the \emph{support} of a vector in tropical space $\RR^d$
 as the set of its zero coordinates.
We define the \emph{support poset} of a matrix $M$ to 
be the set of all unions of supports of column vectors of
$M$. This set is partially ordered by inclusion.

\begin{prop}\label{prop:01rank}
The tropical rank of a zero-one matrix with no column of all ones equals the maximum length of a chain
in its support poset.
\end{prop}

The assumption that there is no column of all ones is needed for the 
statement to hold because a column of zeroes
and a column of ones represent the same point in
tropical projective space $\TP^{d-1}$.

\begin{proof}
There is no loss of generality in assuming that every
union of supports of columns of $M$ is actually the support of a
column. Indeed, the tropical sum of a set of columns gives a column whose
support is the union of supports, and appending this column to $M$
does not change the tropical rank since the tropical convex hull 
of the columns remains the same. Therefore, if there is a chain with 
$r$
elements in the support poset we may assume that there is a set of
$r$ columns with supports contained in one another. Since there
is no column of ones, from this we can easily extract an $r\times r$ minor with
zeroes on and below the diagonal and 1's above the diagonal, which is
tropically non-singular.

Reciprocally, suppose there 
is a tropically non-singular $r\times r$ minor $N$. 
We claim that the support poset of $N$ has a chain of length $r$, from which it follows that the support poset
of $M$ also has a chain of length $r$. Assume without loss of generality that the unique minimum permutation sum is
obtained in the diagonal.  This minimum sum cannot be more than one, because if $n_{ii}$ and $n_{jj}$ are both 1 then
changing them for $n_{ij}$ and $n_{ji}$ does not increase the sum. If the minimum is zero, orienting
an edge from $i$ to $j$ if entry $ij$ of $N$ is zero yields an acyclic digraph, which admits an ordering. Rearranging
the rows and columns according to this ordering yields a matrix with 1's above the diagonal and 0's on and below the 
diagonal.
The tropical sum of the last $i$ columns (which corresponds to union of the corresponding supports) then produces a
vector with 0's exactly in the last $i$ positions. Hence, there is a proper chain of supports of length $r$.

If the minimum permutation sum in $N$ is $1$, then
let $n_{ii}$ be the unique diagonal entry equal to $1$.
The $i$-th row in $N$ must consist of all $1$'s: 
if $n_{ij}$ is zero, then changing $n_{ij}$ and $n_{ji}$ for
$n_{ii}$ and $n_{jj}$ does not increase the sum. 
Changing this row of ones to a row of zeroes
does not affect the support poset of $N$ (it just adds an 
element to every support), and yields a non-singular zero-one matrix
with minimum sum zero to which we can apply the argument
in the previous paragraph.
\end{proof}


\begin{ex}
\rm
\label{ex:trank}
The tropical rank of the classical identity 
matrix $C_n$ equals two (for all $n$), since all of
its $3\times 3$ minors are tropically singular, while 
the principal $2\times 2$ minors are not.
The supports of its columns are all the sets of cardinality
$n-1$ and the support poset consists of them and the whole set
$\{1,\dots,n\}$. The maximal chains in the poset have indeed length two.
\end{ex}

As with the matrices of Barvinok rank $r$, the 
$d\times n$ matrices of tropical rank at most $r$ form a
polyhedral fan given as the intersection of the tropical 
hypersurfaces  $\mathcal{T}(f)$ where $f$ runs over 
the {\bf set of} $(r+1) \times(r+1)$-subdeterminants 
of a $d \times n$-matrix of unknowns $(x_{ij})$. 
Note that this is very similar to the Kapranov rank; 
by Theorem~\ref{thm:kaprlift}, the set of 
$d \times n$ matrices of tropical rank is the
intersection of the tropical 
hypersurfaces  $\mathcal{T}(f)$ where $f$ runs over 
the {\bf ideal generated by the} $(r+1) \times(r+1)$-subdeterminants 
of a $d \times n$-matrix of unknowns $(x_{ij})$. 

However, these are not equal; matrices can have Kapranov rank 
strictly bigger than their tropical rank, 
as will be seen in Section~\ref{sec:matroids}. In this sense, the 
subdeterminants of a given size $r \geq 4$ do not form a 
tropical 
basis for the ideal they generate.

\section{Mixed Subdivisions and Corank One}

A useful tool in tropical convexity is the computation of tropical
convex hulls by means of
 mixed subdivisions of the Minkowski sum of several copies of a simplex. We recall the definition of mixed subdivisions, adapted to 
the case of interest to us. See \cite{S} for more details.

\begin{defn}
Let $\Delta^{d-1}$ be the standard $(d-1)$-simplex in $\RR^d$, with vertex set
$A=\{e_1,\dots,e_d\}$. Let $n\Delta^{d-1}$ denote its dilation by a factor of
$n$, which we regard as the convex hull of the Minkowski sum $A+A+\cdots+A$ ($n$ times).
Let  $M=(v_{ij})\subset\RR^{d\times n}$ be a matrix.
Consider the lifted simplices
\[
P_i \quad := \quad \conv
\bigl\{(e_1,v_{1i}),\dots,(e_d,v_{di}) \bigr\}\,\subset\, \RR^{d+1}
\qquad \hbox{for} \,\, i = 1,2,\ldots,n.
\]
The \emph{regular mixed subdivision} of $n\Delta^{d-1}$ induced by $M$
is the set of projections of the lower faces of
the Minkowski sum $P_1+\cdots +P_n$. Here, a face is called lower if
its outer normal cone contains a vector with last coordinate negative.
\end{defn}

It was shown in~\cite[\S 4]{DS} that there is a bijection between the 
cells $X_S$ in the convex hull of the columns of 
$M$ and the interior cells in the regular
subdivision of a product of simplices induced by $M$. Via the Cayley trick~\cite{S}, the latter biject to
interior cells in the regular mixed subdivision defined above. Here we provide a short direct proof of the composition of
these two bijections:

\begin{lem}
\label{lem:mixedtypes}
Let $M\subset\RR^{d\times n}$ and let
$S=(S_1,\dots,S_d)$, where each $S_j$ is a subset of
$\{1,\dots,n\}$. Then, the following properties are equivalent:
\begin{enumerate}
\item There exists a point in $\RR^d$ of type $S$ relative to the $n$ points given 
by the columns of $M$.
\item There is a non-negative matrix $M'$ such that $M'$ is obtained from $M$ by adding constants to rows or columns 
of $M$, and such that $M'_{ji}=0$ precisely when $i\in S_j$.
\item The regular mixed subdivision of $n\Delta^{d-1}$ induced by $M$
has as a cell the Minkowski sum $\tau_1+\cdots+\tau_n$ where
$\tau_i=\conv(\{e_j: i\in S_j\})$.
\end{enumerate}
Moreover, if this happens, the cells referred to in parts (1) and (3)
have complementary dimensions.
\end{lem}

\begin{proof}
Adding a constant to a row of $M$ amounts to 
translating the set of $n$ points in $\TP^{n-1}$,
while adding a constant to a column leaves the point set
unchanged. Consider a cell $X_S$ in the tropical convex hull, 
let $x$ be any point in the relative interior of $X_S$ 
and let $M'$ be the (unique)
matrix obtained by translating the point set by a vector $-x$ and
normalizing every column by adding a scalar so that its minimum coordinate 
equals $0$. Conversely, for a matrix $M'$ as in (2), consider the point 
$x$ whose coordinates are the amounts added to the columns of 
$M$ to obtain $M'$.  The point $x$ is in the tropical convex hull
of the columns of $M$. Let $S$ be its type.
Then  the modified matrix 
$M'$ has zeroes precisely in entries $(j,i)$ with $i\in 
S_j$, proving the equivalence of (1) and (2). 

For the equivalence of (2) and (3), observe that adding a constant to
a row or column of $M$ does not change the mixed subdivision of $\sum
P_i$. For a non-negative matrix $M'$ with at least a zero in every
column, the positions of the zero entries define the face of $\sum P_i$
in the negative vertical direction. Conversely, for every cell of the
regular mixed subdivision, we can apply a linear transformation changing only the last coordinate to
give that cell height zero and all other vertices positive height
(this is what it means to be in the lower envelope.) The resulting
height function is precisely the matrix $M^\prime$ in (2), which
proves the equivalence of (2) and (3).  The assertion on
dimensions is easy to prove.
\end{proof}

This lemma implies that the tropical convex hull is dual to 
the regular mixed subdivision.

\begin{cor}
\label{cor:mixed}
Given a matrix $M$, the poset of types in the tropical convex hull of
its columns and the poset of interior cells 
of the corresponding regular mixed subdivision are
antiisomorphic. 
\end{cor}

\begin{proof}
>From the proof of Lemma~\ref{lem:mixedtypes}, it is clear that the poset of types (under $S<T$ if 
$S_j\subset T_j$ for each $j$) and the poset of cells in the regular mixed subdivision are 
antiisomorphic. Meanwhile, a type $S$ is in the tropical convex hull of its columns if and only if each 
$S_j$ is nonempty; this is the same condition categorizing when the corresponding cell is contained in 
the boundary of the mixed subdivision (which occurs whenever there exists a vertex appearing in no 
summand.)
\end{proof}

\begin{cor}
\label{cor:duality}
Let $M\subset\RR^{d\times n}$.
The tropical rank of $M$ equals $d$ minus the minimal dimension of
an interior cell in the regular mixed subdivision
of $n\Delta^{d-1}$ induced by $M$.
\end{cor}

We can use these tools to prove that the
tropical and Kapranov ranks of a matrix coincide if the latter is
maximal.

\begin{thm}
\label{thm:maxrank}
If a $d\times n$ matrix $M$ has Kapranov rank equal to $d$, then it has
tropical rank equal to $d$ as well.
\end{thm}

\begin{proof}
By Corollary \ref{cor:duality},
 $M$ has tropical rank $d$  if and only if the corresponding regular mixed subdivision has an
interior vertex. The theorem then follows from the next two lemmas.
\end{proof}

\begin{lem}
\label{lem:maxKapranov}
A $d \times n$-matrix $M$ has Kapranov rank less than $d$
if and only if the corresponding regular mixed subdivision has
a cell that intersects all facets of $n\Delta^{d-1}$.
\end{lem}

\begin{proof}
If $M$ has Kapranov rank less than $d$, then its column vectors
lie in a 
tropical hyperplane.  Since all tropical hyperplanes are translates of
one another, there is no loss of generality in assuming that it is the
hyperplane defined by $x_1\oplus \cdots \oplus x_d$. That is, after
normalization, all columns of $M$ are non-negative and have at least two
zeroes. Then, by Lemma~\ref{lem:mixedtypes}, the zero entries of $M$ define a cell $B$ in the regular mixed subdivision
none of whose Minkowski summands are single vertices. In
particular, for every facet $F$ of $\Delta^{d-1}$ and for every
$i\in\{1,\dots,n\}$, the $i$-th summand of $B$ is at least an edge and
hence it intersects $F$. Hence, $B$ intersects all facets of
$n\Delta^{d-1}$. For the converse suppose
the regular mixed subdivision has a cell $B$ which intersects all facets of
$n\Delta^{d-1}$. We may assume that
$M$ gives height zero to the points in that cell and positive height to all
the others. The intersection of $B$ with the $j$-th facet is given by the 
zero entries in $M$ after deletion of the $j$-th row. In particular, 
$B$ intersects the $j$-th facet if and only if every column has a zero entry 
outside of the $j$-th row, and so $B$ intersects all facets if and only if 
all columns of $M$ have at least two zeroes, implying that these all lie in the hyperplane defined by $x_1\oplus \cdots 
\oplus x_d$.
\end{proof}

The cell in the preceding statement need not be
unique. For example, if a tetrahedron is sliced by planes parallel to
two opposite edges, then each maximal cell  meets all the
facets of the tetrahedron.

\begin{lem}
\label{lem:sperner}
In every polyhedral subdivision of a simplex which has no interior
vertices, but arbitrarily many vertices on the boundary,
 there is a cell that intersects all of the facets.
\end{lem}

\begin{proof}
Observe that there is no loss of generality in assuming that
the polyhedral subdivision $S$ is a
triangulation. For a triangulation, we use Sperner's Lemma \cite{AZ}: ``if the
vertices of a triangulation of $\Delta$ are labeled so that (1) the
vertices of $\Delta$ receive different labels and (2) the vertices in
any face $F$ of $\Delta$ receive labels among those of the vertices of
$F$, then there is a fully labeled simplex''.

Our task is to give our triangulation a Sperner labeling with the
property that every vertex labeled $i$ lies in the $i$-th facet of the
simplex. The way to obtain this is: the vertex opposite to facet $i$
is labeled $i+1$. More generally, the label $i$ of a vertex $v$ is
taken so that $v$ is contained in facet $i$ but not on facet
$i-1$. All labels are modulo $d$.
\end{proof}

\section{Matrices of rank two }
\label{sec:rank2}

By  Theorem \ref{thm:detrank}, if a matrix
has tropical rank two, then the
tropical convex hull of its columns is 
one-dimensional. Since it is contractible \cite{DS}, 
this tropical polytope  is a tree. Another way of
showing this is via the corresponding regular mixed
subdivision. Tropical rank 2 means that all the interior cells have
codimension zero or one. Hence, the subdivision is constructed by slicing
the simplex via a certain number of hyperplanes (which do not meet
inside the simplex) and its dual graph is a tree. 
The special case when the matrix  has
Barvinok rank two is characterized by 
the following proposition.

\begin{prop}
\label{prop:brank2}
The following are equivalent for a matrix $M$:
\begin{enumerate}
\item It has Barvinok rank $2$.
\item All its $3\times 3$ minors  have Barvinok rank $2$.
\item The tropical convex hull of its columns is a path.
\end{enumerate}
\end{prop}

\begin{proof}
(1)$\Rightarrow$(2) is trivial (the Barvinok rank of a
minor cannot exceed that of the whole matrix) and
(3)$\Rightarrow$(1) is easy: if a tropical polytope is a path,
then it is the tropical convex hull of its two endpoints.
Proposition~\ref{prop:barvinok} then implies that
the Barvinok rank is two.

For (2)$\Rightarrow$(3) first observe that the case where $M$ is
$3\times3$ again follows from Proposition~\ref{prop:barvinok}. We next
prove the case where $M$ is $d\times 3$ by contradiction: since the
tropical convex hulls of rows and of columns of a matrix are
isomorphic as cell complexes \cite[Theorem 23]{DS}, assume that the
tropical convex hull of the {\em rows} of $M$ is not a path. Then, there are
three rows whose tropical convex hull is not a path, and their
$3\times 3$ minor has Barvinok rank 3. Finally, if $M$ is of arbitrary
size $d\times n$ and the tropical convex hull of its columns is not a
path, consider three columns whose tropical convex hull is not a path
and apply the previous case to them.
\end{proof}

Our goal in this section is to show that if $M$ has tropical rank $2$ then
it has Kapranov rank $2$. 
Following Theorem \ref{thm:kaprlift} (c), 
this is done by constructing an explicit lift
to a rank $2$ matrix over $\tilde K$.

\begin{lem}
\label{lem:decomposition}
Let $M$ be a matrix of tropical rank two.
Let $x$ be a point in the tropical convex hull of the columns of $M$. Let $M'$ be the matrix obtained by
adding $-x$ to every column and then normalizing columns to have
zero as their minimal entry. After possibly reordering the rows
and columns, $M'$ has the following block structure:
\[
M':=
\left(
\begin{array}{ccccc}
{\bf 0} & {\bf 0} & {\bf 0} & \cdots & {\bf 0} \\
{\bf 0} & A_1 & {\bf 0}  & \cdots & {\bf 0}  \\
{\bf 0} & {\bf 0}  & A_2 & \cdots & {\bf 0}  \\
{\vdots} & \vdots & \vdots &\ddots &\vdots \\
{\bf 0} & {\bf 0}  & {\bf 0}  & \cdots & A_k \\
\end{array}
\right),
\]
where the matrices $A_i$ have all entries positive and every $2\times
2$ minor has the property that the minimum of its four entries is
achieved twice. Each {\bf 0}  represents a matrix of zeroes of the
appropriate size, and the first row and column blocks of $M^\prime$ 
may have size zero.
Moreover, the tropical convex hull of the columns of 
$M^\prime$ is the union of
the tropical convex hulls of the column vectors of the 
blocks augmented by the zero vector ${\bf 0}$, and two of
these $k$ trees meet only at the point ${\bf 0}$.
\end{lem}

\begin{proof}
First, adjoin the column $x$ to our matrix if it does not already exist; since $x$ is in the convex hull of $M$, this 
will not change the tropical convex hull of the columns of $M$. We can then simply remove it at the end, when it is 
transformed into a column of all zeroes. Thus, we can assume that one of the columns of the matrix $M^\prime$ consists of 
all zeroes.

The asserted block decomposition
means that any two given columns of $M'$ have
either equal or disjoint cosupports, where the cosupport of a column is
the set of positions where it does not have a zero. To prove that this
holds, just observe that if it didn't then $M'$ would have the following
minor, where $+$ denotes a strictly positive entry. 
(Recall that each column has a zero in it.)
\[
\left(
\begin{array}{ccc}
0 & + & + \\
0 & 0 & + \\
0 & ? & 0 \\
\end{array}
\right)
\]
But this $3 \times 3$-matrix is tropically non-singular.
The assertion
of the $2\times 2$ minors follows from the fact that the
non-negative matrix
\[
\left(
\begin{array}{ccc}
0 & a & b \\
0 & c & d \\
0 & 0 & 0 \\
\end{array}
\right)
\]
is tropically singular if and only if the minimum of $a$, $b$, $c$ and
$d$ is achieved twice.

Finally, the assertion about the convex hulls is trivial, since any linear combination of column vectors from a given 
block will have all zero entries except in the coordinates corresponding to that block. Any path joining two such 
points from different blocks will pass through the origin.
\end{proof}

We next introduce a technical lemma for making
a power series lifting sufficiently generic.

\begin{lem}
\label{lem:genericlift}
Let $A$ be a non-negative matrix with no zero column
and suppose that the smallest entry in $A$ occurs
most frequently in the first column.
Let $\tilde A$ be the matrix
\[
\left(
\begin{array}{ccc}
0 & {\bf 0}  \\
{\bf 0}  & A \\
\end{array}
\right)
\]
obtained by adjoining a row and a column of zeroes.
If $\tilde A$ has Kapranov rank two, then $\tilde A$ has a rank-2 lift 
$F \in \tilde K^{d \times n}$ in which every
$2\times 2$ minor is non-singular and the $i$-th
column can be written as a linear combination $\lambda_i u_1 + \mu_i u_2$ 
of the first two columns $u_1$ and $u_2$, with $\deg(\lambda_i)\ge\deg(\mu_i)= 0$.
\end{lem}

\begin{proof}
Starting with an arbitrary rank-2 lift $\tilde F$ of $\tilde A$, let $F$ be obtained by
adding to every column a $\tilde K$-linear
combination of the first column of $\tilde F$ with coefficients of
sufficiently high degree (so as to not change the degrees of the entries) but otherwise generic. 
This preserves the degree of every entry and thus the rank of the lift,
but makes every $2\times 2$ minor of $F$ non-singular; by ``generic,'' all we require is that the 
ratio between the coefficients of two columns is not equal to the ratio between those two columns if 
they are scalar multiples of each other. No column of $\tilde F$ is a scalar multiple of its first 
column since no column of $\tilde A$ aside from the first is constant, so no column of $F$ is a scalar 
multiple of the first column either.

Since the lift has rank two and the first two columns are linearly independent, the $i$-th column of 
$F$ is now a $\tilde K$-linear combination
$\lambda_i u_1 + \mu_i u_2$ of the first two columns.
If the degrees of $\lambda_i$ and $\mu_i$ are
different, then their minimum must be zero in order to get a
degree zero element in the first entry of column $i$. 
But then $\deg(\mu_i)>\deg(\lambda_i)=0$ is impossible,
because it would make the $i$-th column of $A$ all zero. Hence
$\,\deg(\lambda_i)>\deg(\mu_i)=0$.

If the degrees are equal, then they are non-positive in order to get
degree zero for the first entry in $\lambda_i u_1 + \mu_i u_2$. But
they cannot be equal and negative, or otherwise entries of positive
degree in $u_2$ would produce entries of negative degree in
$u_i$. Hence, $\deg(\lambda_i)=\deg(\mu_i)=0$ in this case.
\end{proof}

\begin{cor}
\label{cor:glue}
Let $A$ and $B$ be non-negative matrices. Assume that the two matrices
\[
\tilde A:=\left(
\begin{array}{cc}
A & {\bf 0}  \\
{\bf 0}  & 0 \\
\end{array}
\right)
\text{ and }
\tilde B:=\left(
\begin{array}{cc}
0 & {\bf 0}   \\
{\bf 0}  & B \\
\end{array}
\right)
\]
have Kapranov rank equal to 2. Then, the matrix
\[
M:=\left(
\begin{array}{ccc}
A   & {\bf 0} & {\bf 0}   \\
{\bf 0} & { 0}  & {\bf 0}  \\
{\bf 0}  & {\bf 0}   & B \\
\end{array}
\right)
\]
has Kapranov rank equal to 2 as well.
\end{cor}

\begin{proof}
We may assume that neither $A$ nor $B$ has a zero column.
Hence Lemma \ref{lem:genericlift}
applies to both of them.
We number the rows of $M$ from $-k$ to $k'$ and its columns from $-l$ to
$l'$, where $k\times l$ and $k'\times l'$ are the dimensions of
$A$ and $B$ respectively. In this way, $A$ (respectively $B$) is the minor of negative
(respectively, positive) indices. The row and column indexed zero consist of all zeroes.
 To further exhibit the symmetry between $A$ and $B$ the columns and rows in $\tilde A$ will be referred to ``in
reverse''. That is to say, the first and second columns of it are the ones indexed 0 and $-1$ in $M$.

We now construct a lifting $F=(a_{i,j})\in \CC\{\{t\}\}^{d\times n}$ of $M$.
We assume that we are given rank-2 lifts of $\tilde A$ and $\tilde B$
which satisfy the conditions of the previous lemma.
Furthermore, we assume that
the lift of the entry $(0,0)$ is the same in both, which can be achieved by
scaling the first row in one of them.

We use exactly those lifts of $\tilde A$ and $\tilde B$ for the upper-left and bottom-right corner minors of $M$.
Our task is to complete that with an entry $a_{i,j}$ for every $i,j$ with $ij<0$, such that
$\deg(a_{i,j})=0$ and the whole matrix still has rank 2. We claim that it suffices to choose
the entry $a_{-1,1}$ of degree zero and sufficiently generic. That this choice fixes the rest of the matrix is easy to see: 
The entry $a_{1,-1}$ is fixed by the fact that the $3\times 3$ minor
\[
\left(
\begin{array}{ccc}
a_{-1,-1}  & a_{-1,0}  & a_{-1,1}  \\
a_{0,-1}  & a_{0,0}  & a_{0,1}  \\
a_{1,-1}  & a_{1,0}  & a_{1,1}  \\
\end{array}
\right)
\]
needs to have rank 2. All other entries $a_{i,-1}$ and $a_{i,1}$ are
fixed by the fact that the entries $a_{i,-1}$, $a_{i,0}$ and $a_{i,1}$
(two of which come from either $\tilde A$ or $\tilde B$) must satisfy
the same dependence as the three columns of the minor above. For each
$j=-l,\dots, -2$ (respectively $j=2,\dots,l'$), let $\lambda_j$ and
$\mu_j$ be the coefficients in the expression of the $j$-th column of
$\tilde A$ (respectively, of $\tilde B$) as $\lambda_j u_0 + \mu_j
u_{-1}$ (respectively, $\lambda_j u_0+ \mu_j u_{1}$). Then, $a_{i,j}=
\lambda_j a_{i,0}+ \mu_j a_{i,-1}$ (respectively, 
$a_{i,j}= \lambda_j a_{i,0}+ \mu_j a_{i,1}$).

What remains to be shown is
that if $a_{-1,1}$ is of degree
zero and sufficiently generic then all the new entries are of degree
zero too. For this, observe that if $j\in\{-l',\dots,2\}$
(resp. $j\in \{2,\dots,l\}$
then $a_{i,j}$ is of degree zero as long as the coefficient of degree
zero in $a_{i,-1}$ (resp. $a_{i,1}$)
are different from the degree zero
coefficients in the quotient $-\lambda_j a_{i,0}/\mu_j$ (here we are
using the assumption that $\deg(\lambda_j)\ge \deg(\mu_j)\ge 0$). In
terms of the choice of $a_{-1,1} $ this
translates to the following determinant having non-zero coefficient
in degree zero:

\[
\left(
\begin{array}{ccc}
 a_{i,-1} & a_{i,0}  & -\lambda_j a_{i,0}/\mu_j\\
 a_{-1,-1} & a_{-1,0}  & a_{-1,1} \\
a_{0,-1}& a_{0,0}  & a_{0,1}\\
\end{array}
\right),
\qquad
(\text{respectively}\quad
\left(
\begin{array}{ccc}
 a_{0,-1} & a_{0,0}  & a_{0,1}\\
 a_{1,-1} & a_{1,0}  & a_{1,1} \\
-\lambda_j a_{i,0}/\mu_j& a_{i,0}  & a_{i,1}\\
\end{array}
\right)
).
\]
That $a_{-1,1}$ and $a_{1,-1}$ sufficiently generic imply non-singularity of these matrices
follows from the fact that
the following $2\times 2$ minors come from the given lifts of $\tilde A$ and $\tilde B$, hence they
are non-singular:
\[
\left(
\begin{array}{ccc}
 a_{i,-1} & a_{i,0}\\
 a_{0,-1} & a_{0,0} \\
\end{array}
\right),
\qquad
\left(
\begin{array}{ccc}
a_{0,0}    & a_{0,1}\\
a_{i,0}   & a_{i,1}\\
\end{array}
\right).
\]
\end{proof}

\begin{thm}
\label{thm:rank2}
Let $M$ be a matrix of tropical rank $2$. Then its Kapranov
rank equals $2$ as well.
\end{thm}

\begin{proof}
The Kapranov rank of $M$ is always at least the tropical rank, so we need only show that the Kapranov rank is less 
than or equal to 2. If the tropical convex hull $P$ of the columns 
of $M$ is a path, then $M$ has Barvinok rank $2$
(by Proposition~\ref{prop:brank2}) and thus  Kapranov rank $2$. 
Otherwise, let $x$ be a node of degree at least three in 
the tree $P$. We apply
the method of Lemma~\ref{lem:decomposition}. 
Since $x$ has degree at least three, it follows that there are at 
least three blocks $A_i$. In particular, $M$ has
at least three columns.
We induct on the number of columns of $M$.
If $M$ has exactly three columns, then each 
block $A_i$ is a single column, and
every row of $M$ has at most one positive entry.
It is easy to construct an explicit lift of rank $2$: 
in each row, lift the positive entry $\alpha$ as
$-t^\alpha$ and the zero entries as $-1$ and $1+t^\alpha$.
If there are rows of zeroes, 
lift them as $(-1,-1,2)$, for example.

Next, suppose that $M$ has $m\ge 4$ columns. The two blocks with the smallest number of
combined columns have at least $2$ and at most $m-2$ rows all together. Possibly after adding a row and column of 
zeroes, this provides a decomposition of our matrix as
\[
M \,\, \, = \,\,\, \left(
\begin{array}{ccc}
{0} & {\bf 0}  & {\bf 0}  \\
 {\bf 0} & A   &{\bf 0}   \\
{\bf 0}  & {\bf 0}   & B \\
\end{array}
\right),
\]
where both $A$ and $B$ have at least two columns ($A$ is the union of these two blocks, $B$ the union 
of all other blocks.) It then follows that the minors
\[
\left(
\begin{array}{ccc}
{0} & {\bf 0}   \\
 {\bf 0} & A   \\
\end{array}
\right)
\qquad
\text{and}
\qquad
\left(
\begin{array}{ccc}
 {0} &{\bf 0}   \\
{\bf 0}     & B \\
\end{array}
\right)
\]
both have fewer columns than the original matrix. By the inductive hypothesis they have Kapranov rank 2. 
Applying Corollary~\ref{cor:glue} completes the inductive step of the theorem.
\end{proof}

In the proof of Lemma~\ref{lem:genericlift} we again
required the ability to pick generic field elements.
Thus, Theorem~\ref{thm:rank2} holds over any infinite coefficient
field, but it may fail over finite fields. This is illustrated by the 
following example. 
Proposition~\ref{prop:TKcompare} and Theorem~\ref{thm:maxrank}  fail here too,
as does the fact that Kapranov rank 
is invariant under insertion of a tropical combination of existing columns. 

\begin{ex}
\label{ex:finitefield}
The matrix
\[
M=\left(
\begin{array}{ccc}
1 & 0 & 0 \\
0 & 1 & 0 \\
0 & 0 & 0 \\
\end{array}
\right)
=
\left(
\begin{array}{ccc}
1 & 0 & 0 \\
2 & 1 & 1 \\
1 & 0 & 0 \\
\end{array}
\right)
\oplus
\left(
\begin{array}{ccc}
1 & 2 & 1 \\
0 & 1 & 0 \\
0 & 1 & 0 \\
\end{array}
\right).
\]
has Barvinok and tropical ranks equal to $2$, but 
Kapranov rank~$3$
over the two-element field~$\FF_2$.
\end{ex}

\section{Matrices constructed from matroids }
\label{sec:matroids}
One of the important properties of rank in usual linear algebra is
that it produces a matroid. Unfortunately, the definitions of tropical
rank, Kapranov rank, and Barvinok rank all fail to do this.  Consider
the configuration of four points in the tropical plane $\TP^2$ given
by the columns of
\[
M=\begin{pmatrix}
0 & 0 & 0 & \phantom{-}0 \\
0 & 0 & 1 & \phantom{-}2 \\
1 & 0 & 0 & -1
\end{pmatrix}.
\]
By any of our three definitions of rank, the maximal independent
sets of columns are $\{1,2\}$, $\{1,3,4\}$, and $\{2,3,4\}$.  These do not all
have the same size, and so they cannot be the bases of a matroid. The
central obstruction here is that the sets $\{1,2,3\}$ and $\{1,2,4\}$ are
(tropically) collinear, but the set $\{1,2,3,4\}$ is not.
Despite this failure, there is a strong connection between tropical 
linear algebra and matroids.

The results in Sections 5 and 6 imply that any matrix whose tropical
and Kapranov ranks disagree
must be at least of size $5\times 5$. The smallest
example we know is  $7\times 7$.  It is based on the 
\emph{Fano matroid}. To explain the example, and to show how to construct many
others, we prove a theorem about tropical representations of
matroids. The reader is referred to~\cite{Ox} for matroid basics.

\begin{defn}
Let $\M$ be a matroid. The {\em cocircuit matrix of $\M$},
denoted $\C(\M)$, has rows indexed by the elements of the ground set
of $\M$ and columns indexed by the cocircuits of $\M$. 
It has a $0$ in entry $(i,j)$ if the $i$-th element is in the $j$-th
cocircuit and a $1$ otherwise.
\end{defn}

In other words, $\C(\M)$ is the zero-one matrix whose columns have the
cocircuits of $\M$ as supports. (As before, the support of a column is its
set of zeroes.) As an example, the matrix $C_n$ of
Section 2 is the cocircuit  matrix of the uniform matroid of rank $2$ with
$n$ elements.  Similarly, the cocircuit matrix of the uniform matroid
$U_{n,r}$ has size $n\times{n\choose r-1}$ and its columns are all the
zero-one vectors with exactly $r-1$ ones. The following results show that
its tropical and Kapranov ranks
equal $r$. The tropical polytopes defined by these matrices are the
\emph{tropical hypersimplices} studied in \cite{J}.

\begin{prop}\label{prop:matroidtrank}
The tropical rank of the cocircuit matrix 
$\C(\M)$ is the rank of the matroid $\M$.
\end{prop}

\begin{proof}
This is a special case of Proposition \ref{prop:01rank} because the
rank of $\M$ is the maximum length of a chain of non-zero covectors,
and the supports of covectors are precisely the unions of supports of
cocircuits.  Note that $\C(\M)$ cannot have a column of ones because
every cocircuit is non-empty.
\end{proof}


\begin{thm}\label{prop:matroidkrank}
If the Kapranov rank of $\C(\M)$ over the ground field $k$ 
is equal to the rank of $\M$, then $\M$ is representable
over $k$. If $k$ is an infinite field, then the converse also holds.
\end{thm}

\begin{proof}
Let $\M$ be a matroid of rank $r$ on $\{1,\ldots,d\}$
which has $n$ cocircuits and suppose that $F \in {\tilde K}^{d \times n}$
is a rank $r$ lift of the cocircuit matrix $\C(\M)$. For
each row $f_i$ of $F$, let $v_i \in k^d$ be the vector of constant terms in $f_i \in {\tilde K}^d$.
We claim that $V = \{v_1,\dots,v_d\}$ is a representation of $\M$. First note that
$V$ has rank at most $r$ since every $\tilde K$-linear relation among the 
vectors $f_i$ translates into a $k$-linear relation among the $v_i$.
Our claim says that $\{i_1,\ldots,i_r\}$ is a basis of $\M$  if and only if $\{v_{i_1},\ldots,
v_{i_r}\}$ is a basis of $V$. Suppose $\{i_1,\ldots,i_r\}$ is a basis of $\M$. 
Then, as in the proof of Proposition~\ref{prop:01rank}, we
can find a square submatrix of $\C(\M)$ using rows $i_1,\ldots,i_r$ 
with $0$'s on and below the diagonal and $1$'s above
it.  This means that the lifted submatrix of constant terms is lower-triangular 
with nonzero entries along the diagonal. It implies
that $v_{i_1},\ldots,v_{i_r}$ are linearly independent, and, since
${\rm rank}(V) \leq r$, they must be a basis. We also conclude
${\rm rank}(V) = r$. If $\{i_1,\ldots,i_r\}$ is not a basis in $\M$, 
there exists a cocircuit containing none of them; this means that
some column of $\C(\M)$ has all $1$'s in rows $i_1,\ldots,i_r$. 
Therefore, $f_{i_1},\ldots,f_{i_r}$ all have zero constant term in that 
coordinate, which means that $v_{i_1},\ldots,v_{i_r}$ are all $0$ in that coordinate. 
Since the cocircuit is not empty, not all vectors
$v_j$ have an entry of 0 in that coordinate, and so $\{v_{i_1},\ldots,v_{i_r}\}$ cannot be a basis.
This shows that $V$ represents $\M$ over $k$, which completes the proof of
the first statement in Theorem \ref{prop:matroidkrank}.

For the second statement, let us assume that
$\M$ has no loops. This is no loss of generality because a loop 
corresponds to a row of $1$'s in $\C(\M)$, which does not increase the
Kapranov rank because every column has at least a zero.
Assume $\M$ is representable over $k$ and
fix a $d \times n$-matrix $A \in k^{d \times n}$ such that
the rows of $A$ represent $\M$ and the sets of non-zero coordinates along
the columns of $A$ are the cocircuits of $\M$.
Suppose $\{1,\ldots,r\}$ is a basis of $\M$ and
let $A'$ be the submatrix of $A$ consisting of the first $r$ rows. Write
\[
 A \quad = \quad \begin{pmatrix} {\bf I}_r \\ C \end{pmatrix} \cdot A' 
\]
where ${\bf I}_r$ is the identity matrix and $C \in k^{(d-r) \times
r}$.  Observe that $A$, hence $C$, cannot have a row of zeroes
(because $\M$ has no loops). 
Since $k$ is an infinite field, there exists a matrix $B'
\in k^{r \times n}$ such that all entries of the $d \times r$-matrix
$\,\begin{pmatrix} {\bf I}_r \\ C \end{pmatrix} \cdot B' \,$ are
non-zero. We now define
$$ F \quad = \quad
\begin{pmatrix} {\bf I}_r \\ C \end{pmatrix} \cdot  (A' + t B')
\quad \in \,\, {\tilde K}^{d \times n}. $$
This matrix has rank $r$ and
$\, {\rm deg}(F) \, = \, \C(\M)$. This completes the proof of Theorem
\ref{prop:matroidkrank}.
\end{proof}

If $k$ is representable over a finite field, its Kapranov rank (with respect to that field) may still 
exceed its tropical rank. It is easy to find examples - for example, the matroid represented by 
$\{(0,1),(1,0),(1,1),(0,0)\}$ over $\FF_2$ will work.

\begin{cor}
\label{cor:nonrepresentable}
Let $\M$ be a matroid which is not representable over a
given  field $k$. Then the Kapranov rank with respect to $k$ of the
tropical matrix $\C(\M)$ exceeds its tropical rank.
\end{cor}

\begin{figure}[htb]
 \includegraphics[scale=.75]{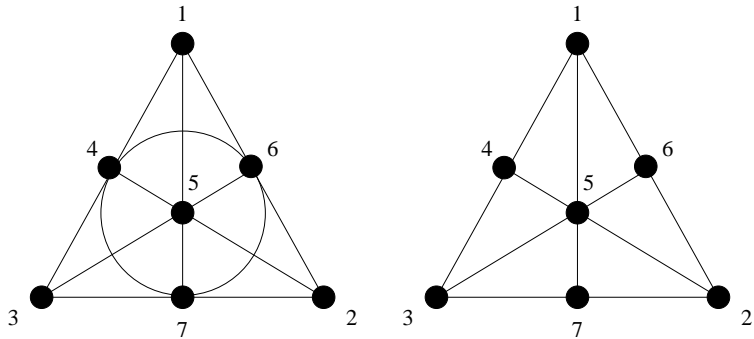}{}
\caption{\label{fig:fano}
The Fano (left) and non-Fano (right) matroids.}
\end{figure}

This corollary furnishes many examples of matrices 
whose Kapranov rank exceeds their tropical rank.
Consider, for example, the Fano and non-Fano matroids,
depicted in Figure~\ref{fig:fano}.
 They both have rank three and seven
elements. The first is only representable over fields of
characteristic two, the second only over fields of characteristic
different from two. In particular, Corollary~\ref{cor:nonrepresentable}
applied to these two matroids implies that over every field
there are matrices with tropical rank equal to three
and Kapranov rank larger than that. Also, it shows that
the Kapranov rank of a matrix may be different over different fields
$k$ and $k'$, even if $k$ and $k'$ are assumed to be 
algebraically closed. This is a more significant discrepancy than that of 
Example~\ref{ex:finitefield}, which used a finite  field.

More explicitly, the cocircuit matrix of the Fano matroid is
\[
\C(\M)=\left(
\begin{array}{ccccccc}
1 & 1 & 0 & 1 & 0 & 0 & 0 \\
0 & 1 & 1 & 0 & 1 & 0 & 0 \\
0 & 0 & 1 & 1 & 0 & 1 & 0 \\
0 & 0 & 0 & 1 & 1 & 0 & 1 \\
1 & 0 & 0 & 0 & 1 & 1 & 0 \\
0 & 1 & 0 & 0 & 0 & 1 & 1 \\
1 & 0 & 1 & 0 & 0 & 0 & 1 \\
\end{array}
\right).
\]
This matrix is the smallest known example of a matrix whose Kapranov
rank over $\CC$ (four) is strictly larger than its tropical
rank (three). Put differently, the seven columns of this matrix (in
$\TP^6$) have as their tropical convex hull a two-dimensional cell
complex which does not lie in any two-dimensional linear subspace of
$\TP^6$, a feature decidedly absent from ordinary
geometry. 

Applied to non-representable matroids, such as 
the \emph{Vamos matroid} (rank 4, 8
elements, 41 cocircuits) or the \emph{non-Pappus matroid} (rank 3, 9
elements, 20 cocircuits) \cite{Ox},
Corollary~\ref{cor:nonrepresentable} yields matrices  
with different Kapranov and
tropical ranks over {\em every} field.
One can also get examples in
which the difference of the two ranks is arbitrarily large. Indeed,
 given matrices $A$ and $B$, we can construct the matrix
\[
M:=
\left(
\begin{array}{cc}
A & \infty \\
\infty' & B \\
\end{array}
\right)
\]
where $\infty$  and $\infty'$
denote matrices of the appropriate dimensions and whose entries are
sufficiently large. 
Appropriate choices of these large values (pick the
extra columns to be points in the tropical convex hull of the columns of $A$ and $B$
and add large constants to 
each column) will 
ensure that the tropical and Kapranov ranks of $M$ are the sums of those of $A$ and of $B$. 
The difference 
between the Kapranov and tropical ranks of $M$ is 
equal to the sum of this difference for $A$ and for $B$. 

The construction in Theorem \ref{prop:matroidkrank} is
closely related to  the  \emph{Bergman complex} of the matroid $\M$.
Ardila and Klivans~\cite{AK} showed that this complex
is triangulated by the order complex of the lattice of flats of $\M$.
Since flats correspond to unions of cocircuits,
the following result is easily derived:

\begin{prop}
The Bergman complex of the matroid $\M$ is equal to the tropical 
convex hull of the rows of the modified 
cocircuit matrix $\C^\prime(\M)$, where the $1$'s in $\C(\M)$ are replaced by $\infty$'s.
\end{prop}

For the Fano matroid, the Bergman complex is the cone over the incidence 
graph of points and lines in the matroid. It consists of $15$
vertices, $35$ edges and $21$ triangles.

\section{Related work and open questions}

As mentioned in the introduction, our definition of non-singular square matrix corresponds to the notion of  ``strongly regular'' in the literature on the max-plus (or min-plus) algebra.
 The definition of ``regular matrix''  in \cite{B,BH,C} is the
following one, for which we prefer to use a different name:

\begin{defn}
A square matrix $M$ is {\em positively tropically regular} 
if, in the formula for its tropical determinant, the minimum over all even permutations equals the minimum over odd permutations. 
The {\em positive tropical rank} of a matrix is the maximum size of a 
positively tropically regular minor.
\end{defn}

The reason for this terminology is that $M$ is positively tropically regular
if it lies outside the positive tropical variety defined by the determinant.
For basics on positive tropical varieties and a detailed study of the
positive tropical Grassmannian see \cite{SW}. The positive
tropicalization of determinantal varieties leads also to a 
notion of \emph{positive Kapranov rank} that satisfies the inequalities
\[
\hbox{pos. tropical rank}\,(M) \quad \leq \quad
\hbox{pos. Kapranov rank}\,(M) \quad \leq \quad
\hbox{Barvinok rank}\,(M) .
\]
Of course, the tropical and Kapranov ranks are less than or equal to their
positive counterparts. 

Our notion of tropical rank, however, appears in \cite{BH,C} under a
different name. Proposition~\ref{prop:SLI} below was 
previously proved in~\cite{BH}:

\begin{defn}
\label{strongregular}
The columns of a matrix $M\in \RR^{d\times n}$ 
are {\em strongly linearly independent}
if there is a column vector $b\in \RR^d$ such that the 
tropical linear system $M\odot x=b$ has a
unique solution $x\in \RR^n$.
A square matrix is {\em strongly regular} if its columns are strongly linearly independent.
\end{defn}

\begin{prop}
\label{prop:SLI}
Strongly regular and tropically non-singular are equivalent,
for a square matrix.
\end{prop}

\begin{proof}
Suppose an $r\times r$ matrix $M$ is tropically non-singular; 
then there is some 
$(r-1)$-dimensional cell $X_S$ in the tropical convex  hull  of  its columns
in $\TP^{r-1}$. After relabeling we have $S_i = \{i\}$ for $i=1,2,\ldots,r$.
 Then taking a point in the relative 
interior of $X_S$ yields
 a vector $b \in \RR^r$ for which $M\odot x = b$ has a unique solution, each $x_i$ being 
necessarily equal to $b_i-m_{ii}$.

Conversely, suppose the columns of an $r \times r$ matrix $M$
 are strongly linearly
independent.  Pick $b \in \RR^r$ such that $M\odot x = b$ has a unique
solution. Then, for each $x_j$, there exists a $b_i$ for
which the expression $\sum M_{ik} x_k$ is uniquely minimized for $k=j$
(otherwise we could increase $x_j$ and get the same value for $M\odot
x)$. This is equivalent to $b$ having type $S$, where $S_j = \{i\}$.
\end{proof}

\begin{cor}
\label{cor:SLI}
The tropical rank of a matrix equals the largest size of a
strongly linearly independent subset of its columns.
\end{cor}

We now discuss some algorithmic issues.  Apart from
Corollary~\ref{cor:SLI}, the main result in \cite{BH} is an $O(n^3)$
algorithm to check strong (i.e., tropical)
regularity of an $n\times n$ matrix.  The
key step is to find a permutation that achieves the minimum in the
determinantal tropical sum, which is the {\em assignment problem} in
combinatorial optimization \cite{PS}. 
Similarly, it is shown in \cite{B} that the
problem of testing positive tropical regularity of square matrices is
equivalent to the problem of testing existence of even cycles 
in directed graphs. 


For the Barvinok rank, we quote some results by \c{C}ela et al.~\cite{CRW}: 

\begin{prop}
The computation of the Barvinok rank of a matrix $M \in \{0,1\}^{d
\times n}$ is an NP-complete problem. Deciding whether a matrix has
Barvinok rank 2 can be done in time $O(dn)$.
\end{prop}

NP-completeness is proved by a reduction to the problem of covering 
a bipartite graph by complete  bipartite subgraphs. For the case of rank 2, 
an algorithm is derived from the fact that matrices of Barvinok rank 2 are
 permuted  \emph{Monge matrices}. 
\c{C}ela et al. also prove 
that a matrix has Barvinok rank $2$ if and only if all its $3\times 3$
minors do (our Proposition \ref{prop:brank2}) and that the Barvinok rank 
is bounded below by the maximum size of a strongly regular minor 
(i.e., by the tropical rank).

\medskip

We finish by listing some open questions, most of them with an
algorithmic flavor:

\begin{enumerate}
\item Singularity of a single minor can be tested in polynomial time.
But a naive algorithm to compute the tropical rank would need to check 
an exponential number of them.
Can the tropical rank of a matrix be computed in polynomial time? 
In other words, is there a tropical analogue of Gauss elimination?

\item Fix an integer $k$. The number of square minors of size at
most $k+1$ of a $d\times n$ matrix $M$ is polynomial in $dn$. Hence,
there is a polynomial time
algorithm for deciding whether $M$ has tropical
rank smaller or equal to $k$. Is the same true for the
Barvinok rank? It is even open whether Barvinok rank equal to 3 can be
tested in polynomial time.

\item For a fixed $k$, a positive answer to either of the following
two questions would imply a positive answer to the previous one:
\begin{enumerate}
\item Is there a number $N(k)$ such that if all minors of $M$ of size
at most $N(k)$ have Barvinok rank at most $k$ then $M$ itself has
Barvinok rank at most $k$?  Proposition~\ref{prop:barvinokCn} shows
that $N(k)\ge{k+1\choose
\left\lfloor\frac{k+1}{2}\right\rfloor}$. 

\item Is there a polynomial time algorithm for the Barvinok rank of matrices with tropical rank bounded by $k$? (This is open even for $k=2$).

\item Can we obtain a bound on the Kapranov rank given the tropical rank? That is, given a positive 
integer $r$, can we find a bound $N(r)$ so that all matrices of tropical rank $r$ have Kapranov rank at 
most $N(r)$? The example of the classical identity matrix shows that the same cannot be done for 
Barvinok rank.

\end{enumerate}

\item Can the Barvinok rank of a matrix $M$ be defined in terms of the
regular mixed subdivision of $n\Delta^{d-1}$ produced by $M$? Ideally,
we would like a ``nice and simple'' characterization such as the one
given for the tropical rank in Corollary~\ref{cor:duality}. But the
question we pose is whether matrices producing the same mixed
subdivision have necessarily the same Barvinok rank.

\item All the questions above are open for the Kapranov rank, too.

\item Does there exist a $5 \times 5$-matrix which has
tropical rank $3$ but Kapranov rank $4$?
\end{enumerate}

\medskip

\noindent
{\bf Acknowledgement:} We thank G\"unter Rote for helpful discussions and for
 pointing  us to references \cite{B,BH}. 


\end{document}